\pdfoutput=1
\RequirePackage{ifpdf}
\ifpdf 
\documentclass[pdftex]{sigma}
\else
\documentclass{sigma}
\fi

\numberwithin{equation}{section}
\numberwithin{theorem}{section}
\numberwithin{proposition}{section}
\numberwithin{lemma}{section}
\numberwithin{corollary}{section}
\numberwithin{definition}{section}
\numberwithin{example}{section}
\numberwithin{remark}{section}
\numberwithin{note}{section}

\usepackage[all]{xy}
\def\C{\mbox{\bbb{C}}}
\def\R{\mbox{\bbb{R}}}
\def\Z{\mbox{\bbb{Z}}}
\def\SZ{\mbox{\sbbb{Z}}}
\def\cd{\C^d}
\def\rd{\R^d}
\def\rn{\R^n}
\def\zd{\Z^d}
\def\td{T^d}
\def\rddu{\big(\rd\big)^*}
\def\rndu{\big(\rn\big)^*}
\def\et1{e^{2\pi i\theta_1}}
\def\etd{e^{2\pi i\theta_d}}
\def\vz{\underline{z}}
\def\zjs{|z_j|^2}
\def\cl{{\cal C}_\lambda}
\def\c0{{\cal C}_0}
\def\D{\Delta}
\def\db{\Delta_{b}^o}
\def\dr{\Delta_{r}^o}
\def\xd{X_1,\ldots,X_d}
\def\ld{\lambda_1,\ldots,\lambda_d}
\def\G{\Gamma}
\def\wt{\tilde{W}}
\def\wsh{W^{\sharp}}

\newfont{\bbb}{msbm10 scaled\magstephalf}
\newfont{\sbbb}{msbm7 scaled\magstephalf}

\begin{document}

\allowdisplaybreaks

\renewcommand{\PaperNumber}{021}

\FirstPageHeading

\ShortArticleName{Ammann Tilings in Symplectic Geometry}

\ArticleName{Ammann Tilings in Symplectic Geometry}

\Author{Fiammetta BATTAGLIA~$^\dag$ and~Elisa PRATO~$^\ddag$}
\AuthorNameForHeading{F.~Battaglia and~E.~Prato}

\Address{$^\dag$~Dipartimento di Matematica e Informatica ``U.~Dini'',
Via S.
Marta 3, 50139 Firenze, Italy}
\EmailD{\href{mailto:fiammetta.battaglia@unifi.it}{fiammetta.battaglia@unifi.it}}
\URLaddressD{\url{http://www.dma.unifi.it/~fiamma/}}

\Address{$^\ddag$~Dipartimento di Matematica e Informatica ``U.~Dini'',
Piazza Ghiberti 27, 50122 Firenze, Italy}
\EmailD{\href{mailto:elisa.prato@unifi.it}{elisa.prato@unifi.it}}
\URLaddressD{\url{http://www.math.unifi.it/people/eprato/}}

\ArticleDates{Received November 09, 2012, in f\/inal form February 27, 2013; Published online March 06, 2013}

\Abstract{In this article we study Ammann tilings from the perspective of symplectic geo\-metry.
Ammann tilings are nonperiodic tilings that are related to quasicrystals with icosahedral symmetry.
We associate to each Ammann tiling two explicitly constructed highly singular symplectic spaces
and~we show that they are dif\/feomorphic but not symplectomorphic.
These spaces inherit from the tiling its very interesting symmetries.}


\Classification{53D20; 52C23}

\section{Introduction}
Our general aim is to study the connection between symplectic geometry and~the nonperiodic tilings
that are related to the geometry of quasicrystals\footnote{Quasicrystals are special materials
having discrete nonperiodic dif\/fraction patterns that were experimentally discovered by Shechtman
et al.~\cite{quasicristalli} in 1982.
They have atomic arrangements with symmetries that are not allowed in ordinary crystals.
For a~comprehensive review of this fascinating subject we refer the reader to the recent book by
Steurer and~Deloudi~\cite{steurer}.}.
Considering these tilings from the symplectic viewpoint provides a~concrete way of obtaining new
examples of highly singular symplectic spaces that are endowed with very rich symmetries.
These examples ef\/fectively contribute to understanding the theoretical aspects of the geometry of
this type of singular spaces.
Moreover, we expect that symplectic geometry may be used to shed light on the study of the tilings.

We started this program with the study of two tilings of the plane: Penrose rhombus~\cite{rhombus}
and~kite and~dart~\cite{kite} tilings.

In this article we focus our attention for the f\/irst time on three-dimensional tilings.
We consider {\em Ammann tilings}, which are the three-dimensional analogues of Penrose rhombus
tilings.
They were introduced by Ammann in the 70's~\cite{senechal-ammann} and~turned out to be related to
quasicrystals with icosahedral symmetry~\cite{steinhardt}.
As we will see later, the third dimension yields an initially unexpected richness and~complexity.

The main idea underlying the connection between symplectic geometry and~tilings is
a~ge\-ne\-ralization~\cite{p} of the Delzant construction~\cite{d}, which we use to associate to each
tile an explicitly constructed symplectic space.
We recall that the Delzant construction associates a~symplectic toric $2n$-manifold to each simple
convex polytope in $\big(\R^n\big)^*$, which is rational with respect to a~lattice and~satisf\/ies an
additional integrality condition.
The problem, in the setting of nonperiodic tilings related to quasicrystals, is that either the
tiles are not individually rational, or they are not simultaneously rational with respect to the
same lattice.
In the generalized construction, however, the lattice is replaced by a~{\em quasilattice}
and~rationality is replaced by the notion of {\em quasirationality}.
The resulting space is a~$2n$-dimensional {\em quasifold}, a~generalization of manifolds
and~orbifolds that was f\/irst introduced by the second-named author in~\cite{p}; the group acting
is no longer a~torus but a~{\em quasitorus}~\cite{p}.

Ammann tilings are made of two kinds of tiles: an oblate rhombohedron and~a~prolate rhombohedron
having same edge lengths.
These rhombohedra, although separately rational, are not simultaneously rational with respect to
the same lattice.
However, the geometry of the tiling ensures that it is possible to choose a~quasilattice $F$ having
the property that each rhombohedron of the tiling is quasirational with respect to $F$
(Proposition~\ref{rotazioni}).
We then apply the genera\-li\-zed Delzant construction simultaneously to each rhombohedron and~we show
that there is one symplectic quasifold, $M_b$, associated to each of the oblate rhombohedra of the
tiling and~one symplectic quasifold, $M_r$, associated to each of the prolate rhombohedra
(Theorem~\ref{uguali}).
Both quasifolds are globally the quotient of a~manifold (the product of three $2$-spheres) modulo
the action of a~discrete group.
There is a~unique dif\/feotype and~two symplectotypes associated to the tiling.
In fact, we show that $M_b$ and~$M_r$ are dif\/feomorphic but not symplectomorphic, consistently
with the fact that the two dif\/ferent types of tiles have dif\/ferent volumes (Theorem~\ref{tipo}).

We remark that quasilattices are the fundamental structure underlying both nonperiodic tilings
related to quasicrystals and~the corresponding symplectic quasifolds.
This is particularly evident for Ammann tilings, and~the related physics of icosahedral
quasicrystals.
The novelty here, with respect to two-dimensional tilings, is that, in this context, there are
actually {\em three} important quasilattices: the quasilattice $F$ above, known as {\em
face-centered lattice}, the {\em body-centered lattice}, $I$, and~the {\em simple icosahedral
lattice}, $P$.
These are the only three quasilattices that have icosahedral symmetry~\cite{rmw}.
The quasilattice $P$, up to a~suitable rescaling, has the property of containing all of the
vertices of the Ammann tiling.
The quasilattices $F$ and~$\frac{1}{2}I$ are usually thought of in physics as the dual of each
other, as they are projections of two $6$-dimensional lattices that are dual to one another in the
standard sense.
Consistently, we show that, in our symplectic setting, $F$ is the group quasilattice of the
quasitorus $D^3=R^3/F$ and that~$\frac{1}{2}I$ can be thought of as its dual, the weight
quasilattice (Section~\ref{geoquasi}).

The paper is structured as follows: in Section~\ref{Delzant} we recall the generalized Delzant
construction; in Section~\ref{quasilattices} we introduce the quasilattices $F$, $I$ and~$P$,
and~we discuss their connection with the tiling; in Section~\ref{interpretation} we construct the
symplectic quasifolds $M_b$ and~$M_r$; in Section~\ref{geoquasi} we study their local geometry;
f\/inally, in Section~\ref{tipi} we show that they are dif\/feomorphic but not symplectomorphic.

\section{The generalized Delzant construction}
\label{Delzant}
We now recall from~\cite{p} the generalized Delzant construction.
For the notion of quasifold, of related geometrical objects and~for a~number of examples we
refer the reader to the original article~\cite{p} and~to
\cite{kite}, where some of the def\/initions were reformulated.

Let us recall what a~{\em simple} convex polytope is.
\begin{definition}[simple polytope] A dimension $n$ convex polytope $\D\subset\rndu$ is said to be
{\em
simple} if there are exactly $n$ edges stemming from each
vertex.\end{definition} Let us next def\/ine the notion of {\em
quasilattice}, introduced in~\cite{mackay}:
\begin{definition}[quasilattice]
Let $E$ be a~real vector space.
A {\em quasilattice} in $E$ is the
$\Z$-span of a~set of $\R$-spanning vectors, $Y_1,\ldots,Y_d$, of
$E$.
\end{definition}
Notice that $\hbox{Span}_{\SZ}\{Y_1,\dots,Y_d\}$ is a~lattice if and
only if it admits a~set of generators which is a~basis of $E$.

Consider now a~dimension $n$ convex polytope $\D\subset\rndu$ having
$d$ facets.
Then there exist elements $\xd$ in $\rn$ and~$\ld$ in
$\R$ such that
\begin{gather}
\label{polydecomp}
\D=\bigcap_{j=1}^d\big\{\mu\in\rndu\,|\,\langle\mu,X_j\rangle\geq\lambda_j\big\}.
\end{gather}
\begin{definition}[quasirational polytope] Let $Q$ be a~quasilattice in
$\rn$.
A convex polytope $\D\subset\rndu$ is said to be {\em
quasirational} with respect to $Q$ if the vectors $\xd$ in~\eqref{polydecomp} can be chosen in $Q$.
\end{definition}
We remark that each polytope in $\rndu$ is quasirational with
respect to some quasilattice~$Q$: just take the quasilattice that is
generated by the elements $\xd$ in~\eqref{polydecomp}.
Notice that
if $\xd$ can be chosen in such a~way that they belong to a~lattice,
then the polytope is rational in the usual sense.
Before we go on to
describing the generalized Delzant construction we recall what a
{\em quasitorus} is.
\begin{definition}[quasitorus] Let $Q\subset\rn$ be a~quasilattice.
We call {\em quasitorus} of dimension $n$ the group
and quasifold $D=\rn/Q$.
\end{definition}
For the def\/inition of Hamiltonian action of a~quasitorus on a
symplectic quasifold we refer the reader to~\cite{p}.

For the purposes of this article we will restrict our attention to
the special case $n=3$.
\begin{theorem}[generalized Delzant construction~\cite{p}]\looseness=-1
Let $Q$ be a~quasilattice in $\R^3$ and~let $\D\subset\big(\R^3\big)^*$ be a
simple convex polytope that is quasirational with respect to $Q$.
Then there exists a~$6$-dimensional compact connected symplectic
quasifold $M$ and~an effective Hamiltonian action of the quasitorus
$D=\R^3/Q$ on $M$ such that the image of the corresponding moment
mapping is~$\D$.
\end{theorem}

\begin{proof} Let us consider the space $\cd$ endowed with the standard
symplectic form
\begin{gather*}
\omega_0=\frac{1}{2\pi i}\sum_{j=1}^d dz_j\wedge d\bar{z}_j
\end{gather*}
and the action of the torus $\td=\rd/\zd$ given by
\begin{alignat*}{5}
 &  \tau \colon \  && \td\times\cd && \longrightarrow \ && \cd & \\
&&&  \big(\big(\et1,\ldots,\etd\big),\vz\big) && \longmapsto \ && \big(\et1z_1,\ldots,\etd z_d\big).&
\end{alignat*}
This is an ef\/fective Hamiltonian action with moment mapping given by
\begin{alignat*}{5}
 & J \colon \  && \cd && \longrightarrow \ && \rddu &\\
 &&& \vz &&\longmapsto \ && \sum_{j=1}^d\zjs e_j^*+\lambda, \qquad
\lambda=\mbox{const}\in\rddu. &
\end{alignat*}
The mapping $J$ is proper and~its image is given by the cone
$\cl=\lambda+\c0$, where $\c0$ denotes the positive orthant of
$\rddu$.
Take now vectors $\xd\in Q$ and~real numbers $\ld$ as in~\eqref{polydecomp}.
Consider the surjective linear mapping
\begin{alignat*}{5}
& \pi \colon \ && \rd && \longrightarrow \ && \R^3 &\\
&&& e_j && \longmapsto \ && X_j. &
\end{alignat*}
Consider the dimension $3$ quasitorus $D=\R^3/Q$.
Then the linear
mapping $\pi$ induces a~quasitorus epimorphism $\Pi \colon \td
\longrightarrow D$.
Def\/ine now $N$ to be the kernel of the mapping
$\Pi$ and~choose $\lambda=\sum\limits_{j=1}^d\lambda_j e_j^*$.
Denote by~$i$ the Lie algebra inclusion $\mbox{Lie}(N)\rightarrow\rd$ and
notice that $\Psi=i^*\circ J$ is a~moment mapping for the induced
action of~$N$ on~$\cd$.
Then the quasitorus $\td/N$ acts in a~Hamiltonian fashion on the compact symplectic quasifold
$M=\Psi^{-1}(0)/N$.
If we identify the quasitori $D$ and~$\td/N$ via
the epimorphism $\Pi$, we get a~Hamiltonian action of the quasitorus~$D$ whose moment mapping~$\Phi$ has image equal to
${(\pi^*)}^{-1}(\cl\cap\ker{i^*})=
{(\pi^*)}^{-1}(\cl\cap\mbox{im}\,\pi^*)={(\pi^*)}^{-1}(\pi^*(\D))$
which is exactly $\D$.
This action is ef\/fective since the level set
$\Psi^{-1}(0)$ contains points of the form $\vz\in\cd$, $z_j\neq0$,
$j=1,\ldots,d$, where the $T^d$-action is free.
Notice f\/inally that
$\dim{M}=2d-2\dim{N}=2d-2(d-3)=6$.
\end{proof}

\begin{remark}
\label{choices}
If we want to apply this construction to any simple convex polytope
in $\big(\R^3\big)^*$, then there are two arbitrary choices involved.
The
f\/irst is the choice of a~quasilattice~$Q$ with respect to which the
polytope is quasirational, and~the second is the choice of vectors~$\xd$ in~$Q$ that are orthogonal to the facets of~$\D$ and
inward-pointing as in~\eqref{polydecomp}.
\end{remark}

\section{Ammann tilings and~quasilattices}
\label{quasilattices}
The purpose of this section is to introduce three quasilattices $P$, $F$ and~$I$, that are relevant
for our construction.

Let $\phi=\frac{1+\sqrt{5}}{2}$ be the golden ratio.
We will be using extensively the following fundamental identity
\begin{gather}
\label{laphi}
\phi=1+\frac{1}{\phi}.
\end{gather}

\begin{figure}[t]
\centering
\begin{minipage}[b]{70mm}
\centering
\includegraphics[width=43mm]{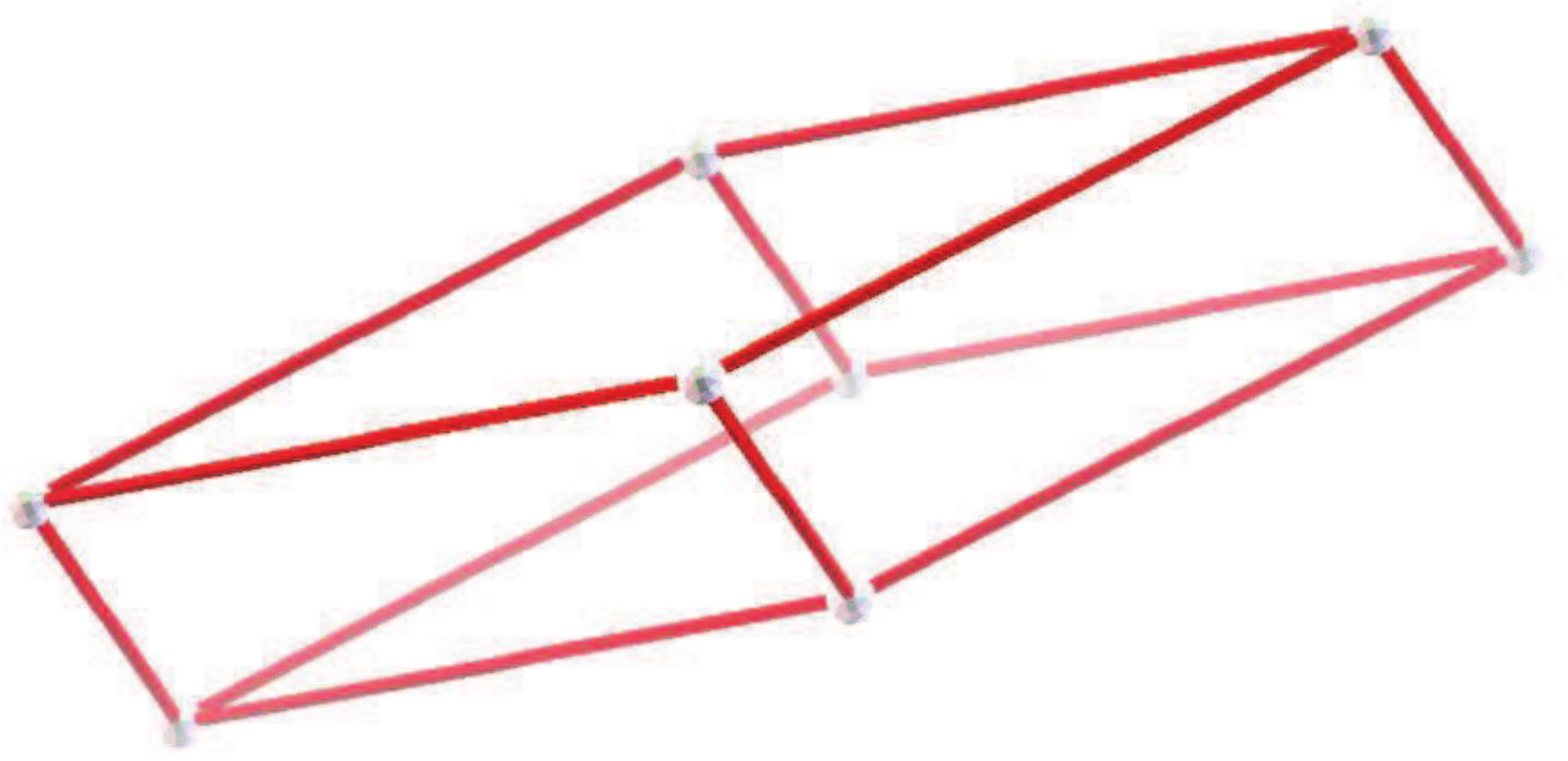}
\caption{The oblate rhombohedron.}\label{oblato}
\end{minipage}
\qquad\quad
\begin{minipage}[b]{70mm}
\centering
\includegraphics[width=43mm]{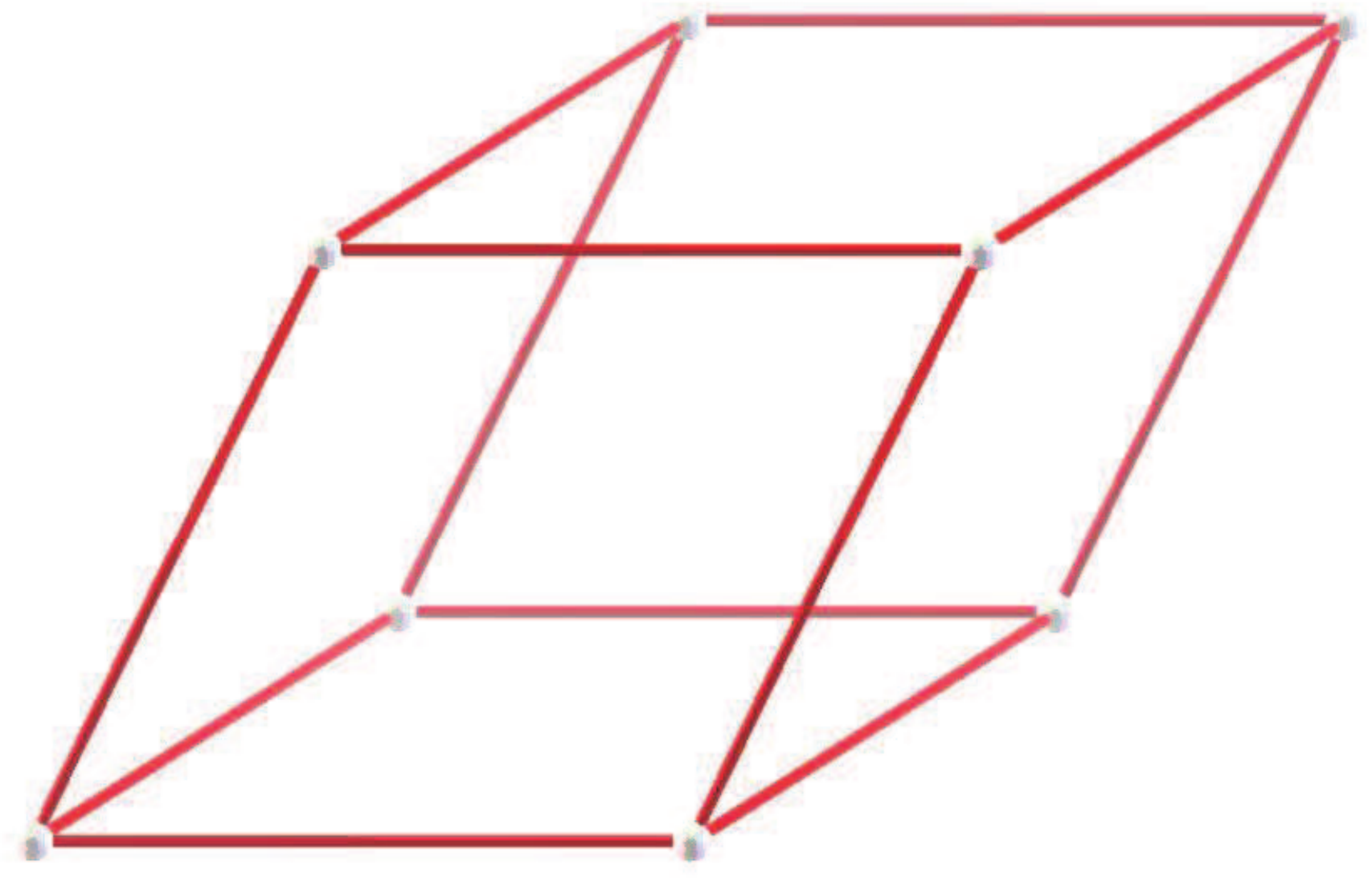}
\caption{The prolate rhombohedron.}\label{prolato}
\end{minipage}
\end{figure}

Let $\sigma$ be a~positive real number and~let us consider an Ammann tiling $\cal T$ with f\/ixed
edge length $\sigma$.
Ammann tilings are nonperiodic tilings of three-dimensional space by
so-called {\em golden rhombohedra}; rhombohedra are called golden
when their facets are given by {\em golden rhombuses}, namely
rhombuses with diagonals that are in the ratio of $\phi$.
There are
two types of such rhombohedra which are called {\em oblate} and~{\em
prolate} (see Figs.~\ref{oblato} and~\ref{prolato})\footnote{All pictures were drawn using the
ZomeCAD software.}.
For a~review of Ammann tilings we refer the reader to \cite{senechal-ammann,steurer}.

Consider now the vectors in $\big(\R^3\big)^*$
\begin{alignat*}{4}
& \alpha_1=\tfrac{1}{\sqrt{2}} (\phi-1,1,0),\qquad  && \alpha_2=\tfrac{1}{\sqrt{2}} (0,\phi-1,1), \qquad && \alpha_3=\tfrac{1}{\sqrt{2}} (1,0,\phi-1),&\\
& \alpha_4=\tfrac{1}{\sqrt{2}} (1-\phi,1,0), \qquad && \alpha_5=\tfrac{1}{\sqrt{2}} (0,1-\phi,1), \qquad && \alpha_6=\tfrac{1}{\sqrt{2}} (1,0,1-\phi). &
\end{alignat*}

\begin{figure}[t]\centering
\begin{minipage}[b]{70mm}
\centering
\includegraphics[width=40mm]{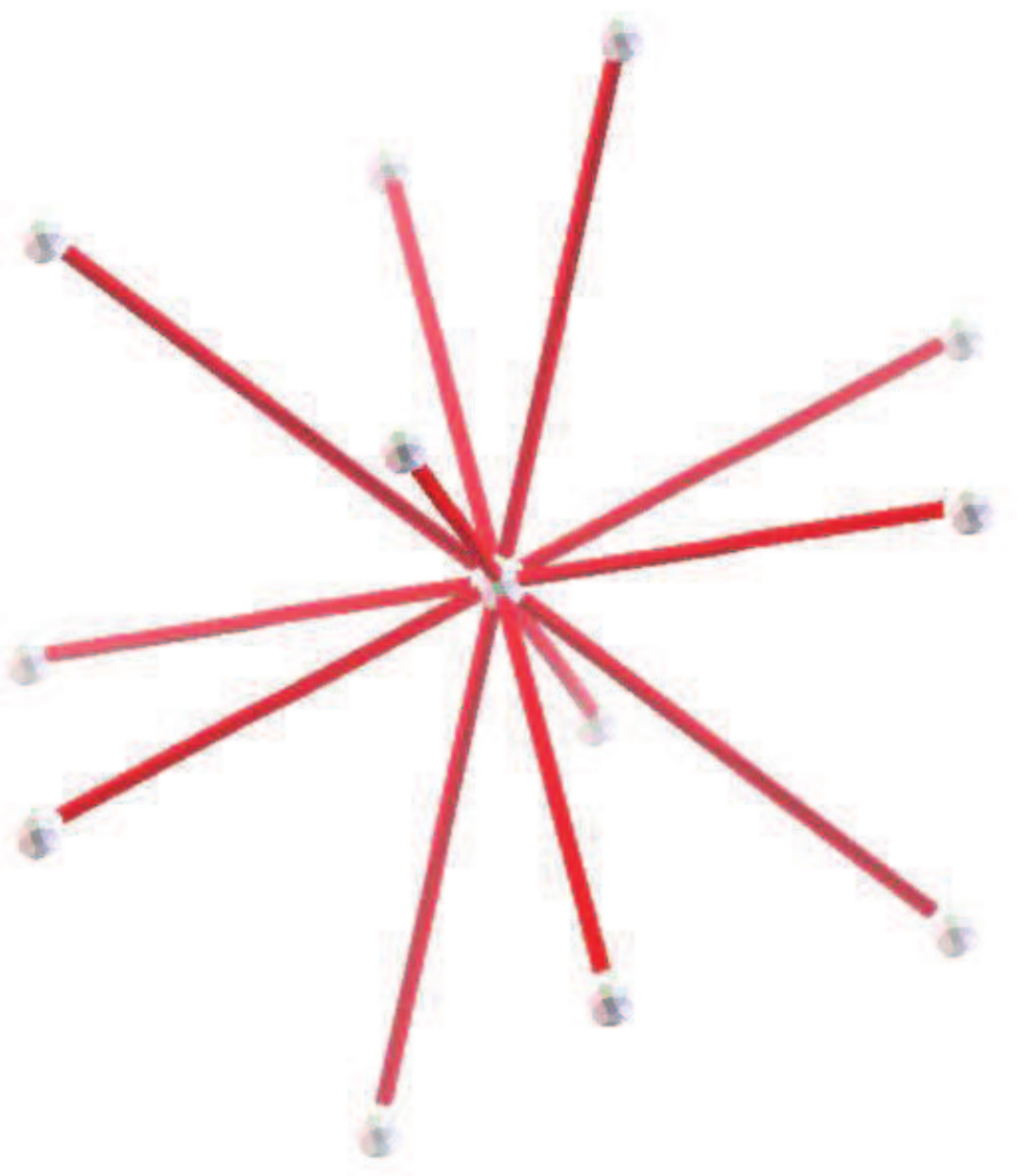}
\caption{The vectors $\pm\alpha_1,\ldots,\pm\alpha_6$.}
\label{stelladirossi}
\end{minipage}
\qquad\quad
\begin{minipage}[b]{70mm}
\centering
\includegraphics[width=45mm]{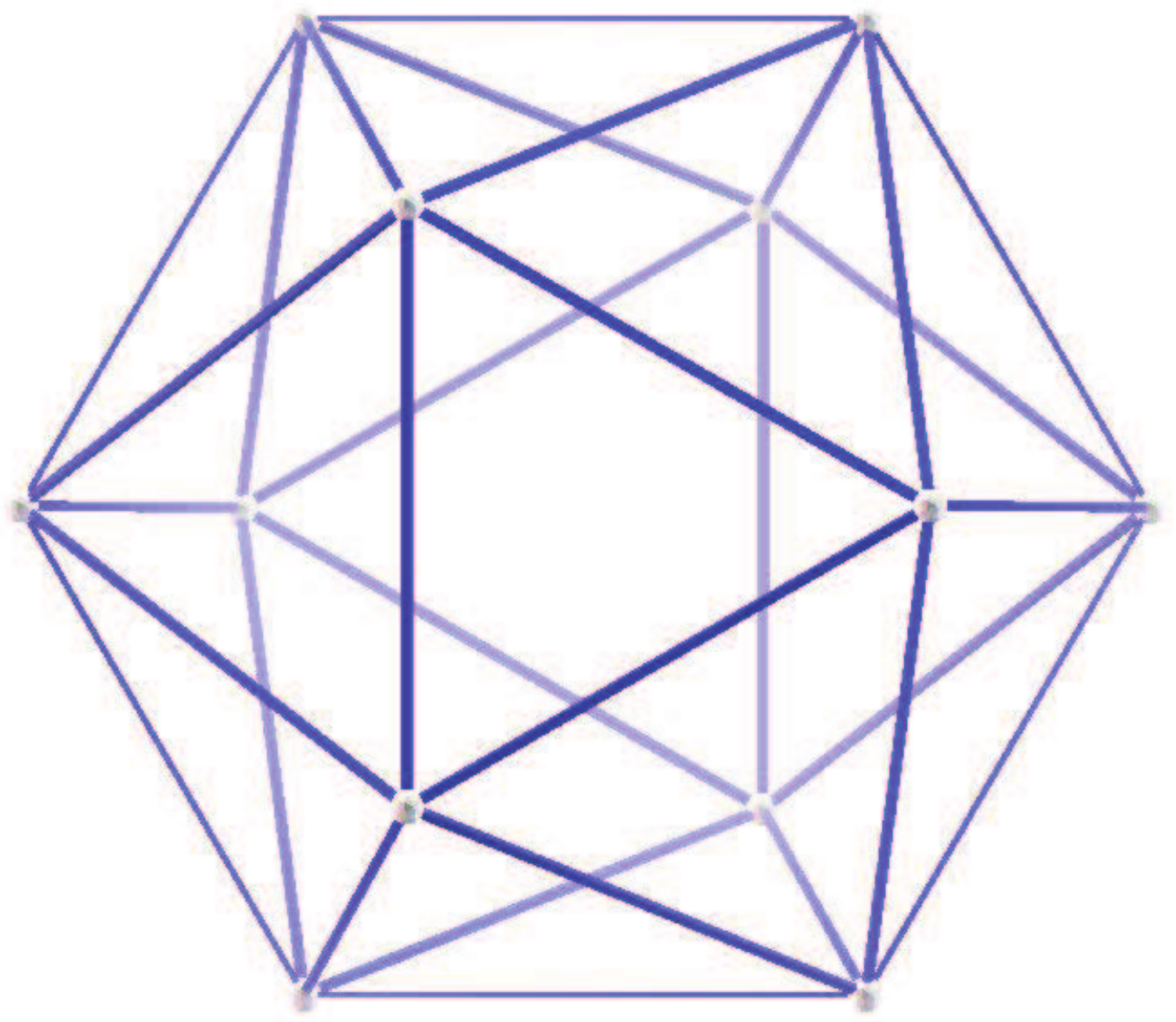}
\caption{The icosahedron.}
\label{icosahedron}
\end{minipage}
\end{figure}

These six vectors and~their opposites point to the twelve vertices
of an icosahedron that is inscribed in the sphere of radius $\sqrt{\frac{3-\phi}{2}}$
(see Figs.~\ref{stelladirossi} and~\ref{icosahedron}); they
generate a~quasilattice $P$ that is known in physics as the {\em simple icosahedral
lattice}~\cite{rmw}.

Let $\delta=\sqrt{\frac{2}{3-\phi}}\sigma$ and
consider the two following golden rhombohedra: the oblate
rhombohedron $\db$, having nonparallel edges
$\delta\alpha_4$, $\delta\alpha_5$, $\delta\alpha_6$, and~the prolate rhombohedron~$\dr$,
having nonparallel edges $\delta\alpha_1$, $\delta\alpha_2$, $\delta\alpha_3$.

Denote by $AB$ one edge of the tiling $\cal
T$.
From now on we will choose our coordinates so that $A=O$ and~so
that $B-A$ is parallel to $\alpha_1$ with the same orientation.
\begin{proposition}
\label{rotazioni}
Let $\cal T$ be an Ammann tiling
with edges of length $\sigma$.
Each vertex of the tiling lies in the
quasilattice $\delta P$.
Moreover, for each oblate rhombohedron $\Delta_b$
in $\cal T$ $($respectively prolate rhombohedron $\Delta_r$ in $\cal
T)$ there is a~rigid motion $\rho$, given by the composition of a
translation with a~transformation of the icosahedral group, such
that $\rho(\Delta_b)$ is $\db$ $($respectively $\rho(\Delta_r)$ is
$\dr)$.\end{proposition}

\begin{proof} Consider f\/irst the icosahedron with its
twenty pairwise parallel facets.
To each pair of parallel facets
there correspond two oblate rhombohedra, one the translate of the
other, and~two prolate rhombohedra, also one the translate of the
other.
Pick one representative for each such couple.
This gives a
total of ten oblate rhombohedra and~ten prolate rhombohedra.
Each of
the ten oblate rhombohedra can be mapped to $\db$ via a
transformation of the icosahedral group, and~in the same way each of
the ten prolate rhombohedra can be mapped to $\dr$.

Now, let $C$ be a~vertex of the tiling that is dif\/ferent from $0$
and the above vertex $B$.
We can join $B$ to $C$ with a~broken line
made of subsequent edges of the tiling.
We denote the vertices of
the broken line thus obtained by
$T_0=A,T_1=B,\dots,T_j,\dots,T_m=C$.
Since the tiles are oblate and
prolate rhombohedra, each vector $Y_j=T_j-T_{j-1}$ is one of the
vectors $\pm\alpha_k$, $k=1,\ldots,6$.
Therefore we have that
$C-A=T_m-T_0=Y_m+\cdots+Y_1$.
This implies that the vertex $C$ lies
in~$\delta P$, that each oblate rhombohedron having $C$ as vertex is the
translate of one of the ten oblate rhombohedra described above and
that each prolate rhombohedron having $C$ as vertex is the translate
of one of the ten prolate rhombohedra described above.
We can
therefore conclude that, for each oblate rhombohedron $\D_b$ having
$C$ as vertex, there exists a~rigid motion $\rho$, given by the
composition of a~translation with a~transformation of the
icosahedral group, such that $\rho(\D_b)=\db$.
The same is true for
the prolate rhombohedra.
\end{proof}

We introduce now a~quasilattice $F\subset\R^3$ with respect to which {\em all} of
the rhombohedra of the tiling are quasirational (cf.\
Remark~\ref{choices}).
This is necessary in order to apply the generalized Delzant
procedure simultaneously to all of the rhombohedra in the tiling.
We take~$F$ to be the quasilattice that is generated by the six vectors{\samepage
\begin{alignat*}{4}
& U_1=\tfrac{1}{\sqrt{2}} (1,\phi-1,\phi),\qquad &&
U_2=\tfrac{1}{\sqrt{2}} (\phi,1,\phi-1), \qquad &&
U_3=\tfrac{1}{\sqrt{2}} (\phi-1,\phi,1), &\\
&  U_4=\tfrac{1}{\sqrt{2}} (-1,\phi-1,\phi),\qquad &&
 U_5=\tfrac{1}{\sqrt{2}} (\phi,-1,\phi-1),\qquad &&
 U_6=\tfrac{1}{\sqrt{2}} (\phi-1,\phi,-1).&
\end{alignat*}
The quasilattice $F$ is known in physics as the {\em face-centered lattice}~\cite{rmw}.}

\begin{figure}[t]\centering
\begin{minipage}[t]{70mm}
\centering \includegraphics[width=45mm]{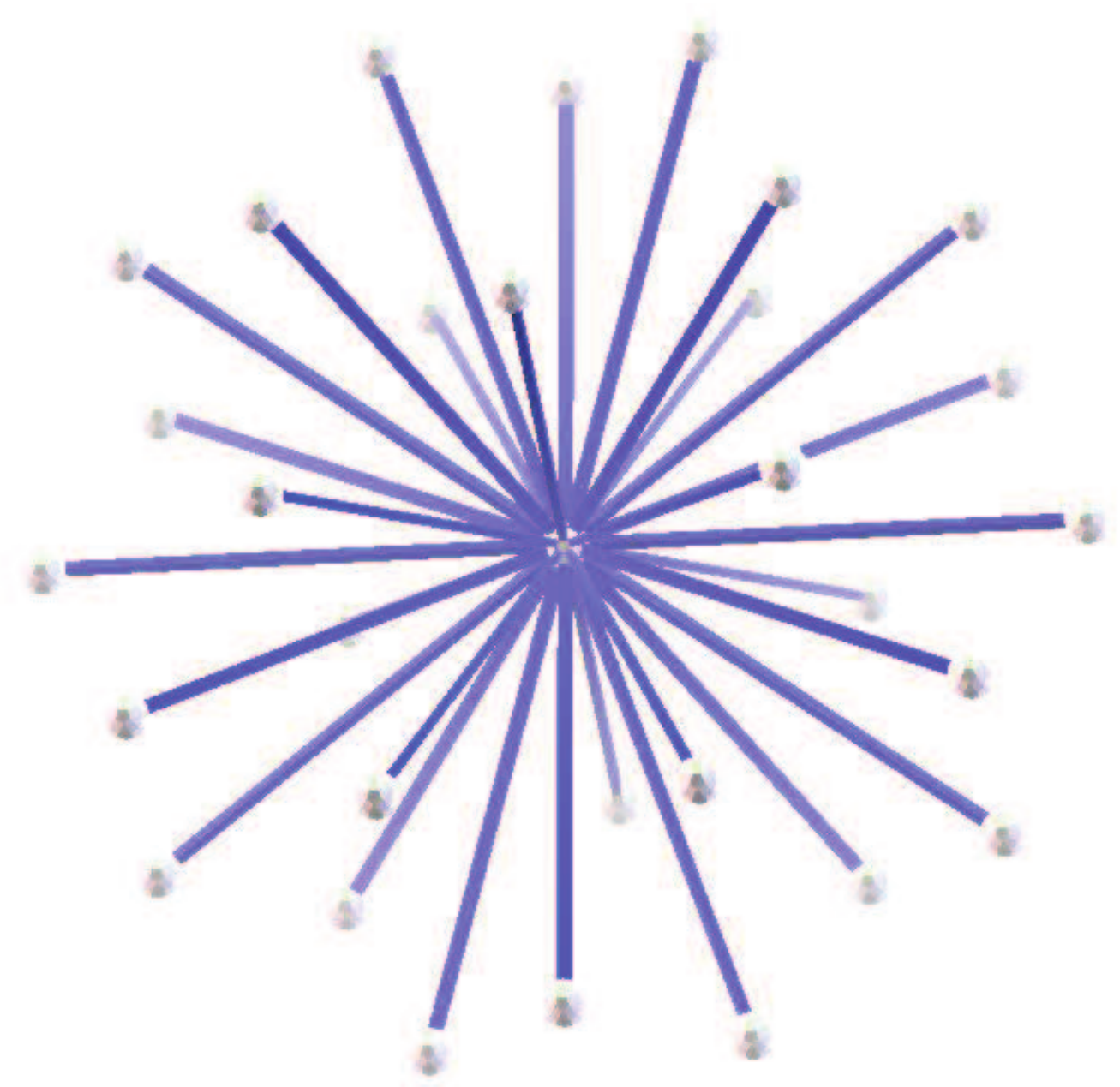}
\caption{The star of norm $\sqrt{2}$ vectors of the quasilattice~$F$.}\label{stelladiblu}
\end{minipage}
\qquad \quad
\begin{minipage}[t]{70mm}
\centering \includegraphics[width=45mm]{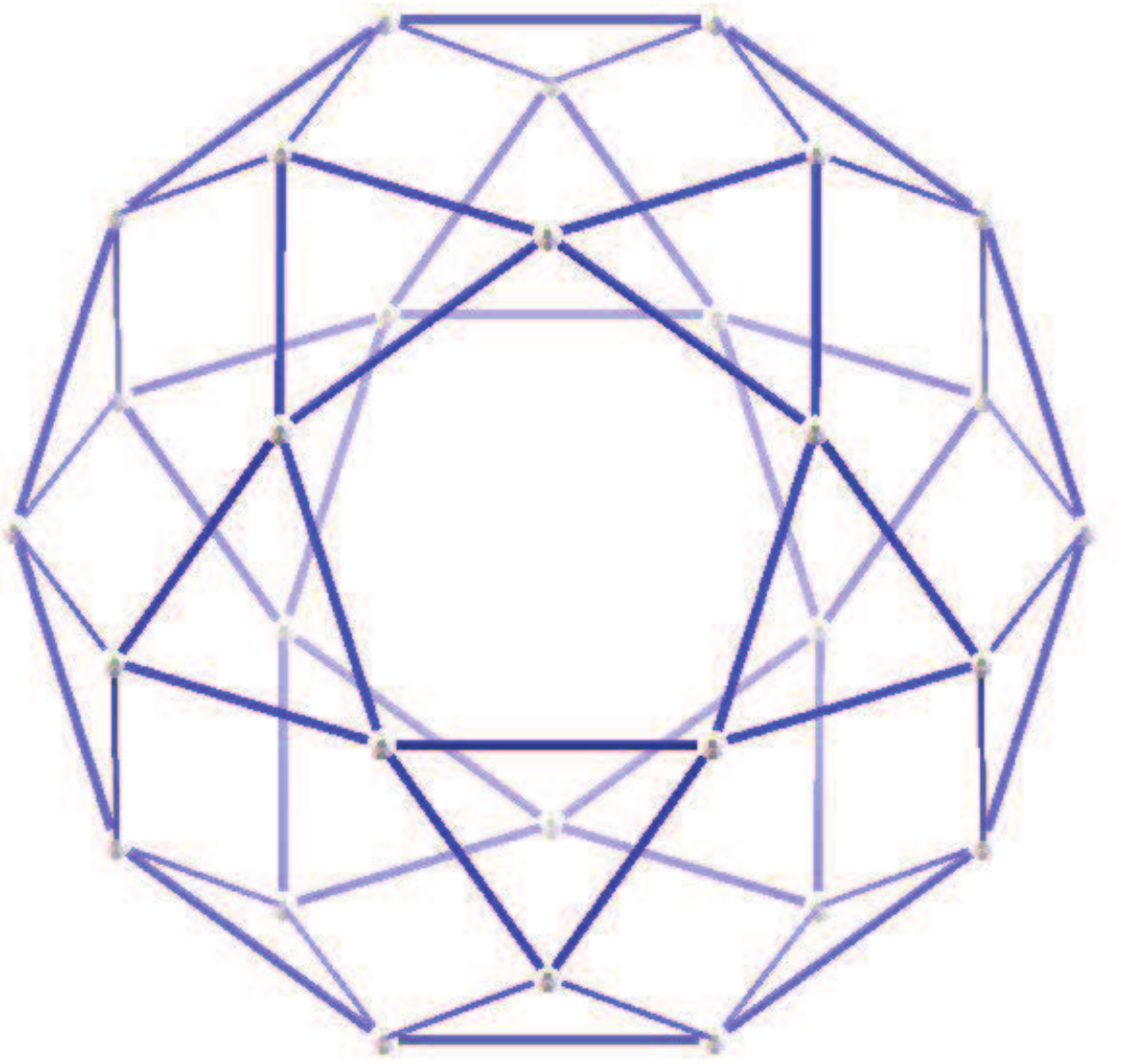}
\caption{The icosidodecahedron.}\label{icosidodecaedro}
\end{minipage}
\end{figure}

The vectors $U_i$ have norm equal to~$\sqrt{2}$.
It can be easily seen that there are exactly $30$ vectors in $F$ having the same norm.
These thirty vectors point to the vertices of an icosidodecahedron inscribed in the sphere of
radius $\sqrt{2}$ (see Figs.~\ref{stelladiblu} and \ref{icosidodecaedro}).

\begin{remark}
Proposition~\ref{rotazioni} implies that, for each facet of the
Ammann tiling, there is a~pair of vectors $\{\alpha_i,\alpha_j\}$
such that the given facet is parallel to the plane $\Pi_{ij}$ generated by
$\{\alpha_i,\alpha_j\}$.
We have $15$ such possible pairs
$\{\alpha_i,\alpha_j\}$, with $i,j=1,\ldots,6$, $i\neq j$.
For each one of them,
two of the $30$ vectors above are orthogonal to the corresponding plane~$\Pi_{ij}$.
This ensures that all of the rhombohedra of the tiling are quasirational with respect to~$F$.
\end{remark}

Another quasilattice that will be useful in the sequel is the quasilattice $I\subset \big(\R^3\big)^*$ that
is generated by
the vectors $\{2\alpha_1,2\alpha_2,2\alpha_3,2\phi\alpha_4,2\phi\alpha_5,2\phi\alpha_6\}$.
The quasilattice $I$ is known in physics as the
{\em body centered lattice}~\cite{rmw}.
\begin{remark}
\label{duali}
The quasilattices $P$, $F$ and~$I$ are invariant under icosahedral
symmetries and~are dense in their respective ambient spaces.
One can show that, if we identify $\R^3$ with its
dual using the standard inner product, we have the following proper inclusions:
\begin{gather}
I\subset F\subset P\subset\tfrac{1}{2}I.
\label{inclusioni}
\end{gather}
Using the notation of Conway--Sloane~\cite{cs}, the lattices $P$, $F$ and~$\frac{1}{2}I$ can be
obtained as the respective projections from the following lattices in $\R^6$: $\Z^6$,
\begin{gather*}
D_6=\hbox{Span}\left\{\sum_{j=1}^6n_je_j\;\Big|\;\sum_{j=1}^6n_j\;\hbox{is even}\right\},
\end{gather*}
and
\begin{gather*}
D_6^*=\hbox{Span}\left\{\frac12\sum_{j=1}^6n_je_j\;\Big|\;n_j\equiv n_k \ (\hbox{mod}\,2)\right\},
\end{gather*}
the projection being given by
\begin{alignat*}{4}
& {\R^6} && \longrightarrow \ && \R^3&  \\
& {e_j} && \longmapsto \ && \alpha_j.&
\end{alignat*}
Coherently with~\eqref{inclusioni} we have the following proper inclusions:
\begin{gather*}
2D_6^*\subset D_6\subset\Z^6\subset D_6^*.
\end{gather*}
The lattice $\Z^6$ is self dual, whilst $D_6$ and~$D_6^*$ are the dual of one another.
A notion of duality
for the icosahedral quasilattices in dimension $3$ is derived
from the above relations of duality in~$\R^6$.
This is coherent with the symplectic setup.
In fact, we will see that the quasilattice $F$ is the group quasilattice
(see Remark~\ref{grouplattice}) and~that the quasilattice~$\frac{1}{2}I$ plays the role of its
dual, the weight quasilattice (see Section~\ref{geoquasi}).
\end{remark}

\section{The tiling from a~symplectic viewpoint}
\label{interpretation}
In this section we perform the Delzant construction to obtain
symplectic quasifolds that can be associated to the oblate and
prolate rhombohedra of an Ammann tiling having edge length $\sigma$.

Let us consider the quasilattice $F$ that we introduced in
Section~\ref{quasilattices}.
As we have seen, all of the rhombohedra
of our tiling are quasirational with respect to~$F$.

\looseness=-1
We begin by considering the oblate rhombohedron $\db$ which has one
of its vertices at the origin and~is determined by the three
non-parallel vectors $\delta\alpha_4$, $\delta\alpha_5$, $\delta\alpha_6$.
This simple
polytope has $6$ facets.
For our construction we choose the $6$
vectors given by $X_1=U_1$, $X_2=U_2$, $X_3=U_3$, $X_4=-U_1$,
$X_5=-U_2$ and~$X_6=-U_3$.
Then the corresponding coef\/f\/icients are
given by $\lambda_1=\lambda_2=\lambda_3=0$ and
$\lambda_4=\lambda_5=\lambda_6=-\frac{\delta}{2\phi}$.
Take now the
surjective linear mapping def\/ined by
\begin{alignat*}{5}
& \pi \colon \ && \R^6 && \rightarrow \ && \R^3 &\\
&&& e_i && \mapsto \ && X_i. &
\end{alignat*}
Its kernel, $\mathfrak n$, is the $3$-dimensional subspace of $\R^6$ that is
spanned by $e_1+e_4$, $e_2+e_5$ and~$e_3+e_6$.
It is the Lie algebra of
$N=\{ \exp(X)\in T^6\,|\, X\in\R^6,\pi(X)\in F \}$.
If $\Psi_b$
is the moment mapping of the induced $N$-action, then
\begin{alignat*}{5}
& \Psi_b \colon \ && \C^6 && \longrightarrow \  && \big(\R^3\big)^* &\\
&&&
\vz && \longmapsto \ && \left(|z_1|^2+|z_4|^2-\tfrac{\delta}{2\phi},
|z_2|^2+|z_5|^2-\tfrac{\delta}{2\phi},|z_3|^2+|z_6|^2-\tfrac{\delta}{2\phi}\right). &
\end{alignat*}
Therefore $\Psi_b^{-1}(0)=S^3_b\times S^3_b\times S^3_b$, where
$S^3_b$ is the sphere in $\R^4$ centered at the origin with radius
$b=\sqrt{\frac{\delta}{2\phi}}$.
In order to compute the group~$N$ we
need the following linear relations between the generators of the
quasilattice $F$:
\begin{gather*}
\left(
\begin{matrix}
U_4
\\
U_5
\\
U_6
\end{matrix}
\right)=\left(\begin{matrix}
1&-\phi&1
\\
1&1&-\phi
\\
-\phi&1&1
\end{matrix}
\right)\left(\begin{matrix}U_1
\\
U_2
\\
U_3\end{matrix}\right).
\end{gather*}
Then a~straightforward computation gives that
\begin{gather*}
N=\left\{ \exp(X)\in T^6\,|\,X= (r+\phi h,s+\phi k,t+\phi l,
r,s,t ),\, r,s,t\in\R,\, h,k,l\in\Z \right\}.
\end{gather*}
We can think of
\begin{gather}
\label{esse1peresse1}
S^1\times S^1\times S^1=\big\{ \exp(X)\in
T^6\,|\,X=(r,s,t,r,s,t),\, r,s,t\in\R \big\}
\end{gather}
as being naturally embedded in~$N$.
The quotient group
\begin{gather*}
\Gamma=\frac{N}{S^1\times S^1\times S^1}
\end{gather*}
is discrete.
In conclusion, the symplectic quotient $M_b$ is given by
\begin{gather*}
M_b=\frac{\Psi_b^{-1}(0)}{N}=\frac{S^3_b\times S^3_b\times S^3_b}{N}=
\frac{S^2_b\times S^2_b\times S^2_b}{\Gamma},
\end{gather*}
where $S^2_b$ is the
sphere in $\R^3$ centered at the origin with radius~$b$.
The
quasitorus $D^3=\R^3/F$ acts on $M_b$ in a~Hamiltonian fashion, with
image of the corresponding moment mapping given exactly by the
oblate rhombohedron~$\db$.

Consider now the prolate rhombohedron $\dr$ that has one vertex in
the origin and~is determined by the three nonparallel vectors
$\delta\alpha_1$, $\delta\alpha_2$, $\delta\alpha_3$.
We now choose the vectors given by
$X_1=U_4$, $X_2=U_5$, $X_3=U_6$, $X_4=-U_4$, $X_5=-U_5$ and
$X_6=-U_6$.
Then the corresponding coef\/f\/icients are given by
$\lambda_1=\lambda_2=\lambda_3=0$ and
$\lambda_4=\lambda_5=\lambda_6=-\frac{\delta}{2}$.
It is immediate to check
that we obtain the same Lie algebra $\mathfrak n$ as in the case of the
oblate rhombohedron.
In order to see what happens to the
corresponding group we need here the inverse relations:
\begin{gather*}
\left(
\begin{matrix}
U_1
\\
U_2
\\
U_3
\end{matrix}
\right)=\left(\begin{matrix}
1&1&\frac{1}{\phi}
\\
\frac{1}{\phi}&1&1
\\
1&\frac{1}{\phi}&1
\end{matrix}
\right)\left(\begin{matrix} U_4
\\
U_5
\\
U_6\end{matrix}\right).
\end{gather*}
To write the relations in this form we used the fundamental identity~\eqref{laphi}.
This identity also implies that we obtain the same
group $N$ as in the case of the oblate rhombohedron.

The moment mapping $\Psi_r$ is given by
\begin{alignat*}{5}
& \Psi_r \colon\ &&  \C^6 && \longrightarrow \ && \big(\R^3\big)^* &\\
&&&
\vz && \longmapsto \ && \left(|z_1|^2+|z_4|^2-\tfrac{\delta}{2},
|z_2|^2+|z_5|^2-\tfrac{\delta}{2},|z_3|^2+|z_6|^2-\tfrac{\delta}{2}\right). &
\end{alignat*}
Therefore
\begin{gather*}
M_r=\frac{\Psi_r^{-1}(0)}{N}=\frac{S^3_r\times S^3_r\times S^3_r}{N}
=\frac{S^2_r\times S^2_r\times S^2_r}{\Gamma},
\end{gather*}
where
$S^2_r\subset\R^3$ and~$S^3_r\subset\R^4$ are the spheres centered
at the origin with radius $r=\sqrt{\frac{\delta}{2}}$.
The quasifold~$M_r$ is acted on by the same quasitorus $D^3=\R^3/F$ that we
obtained for the oblate rhombohedron.
This action is Hamiltonian and
the image of the corresponding moment mapping is exactly the prolate
rhombohedron~$\dr$.

\begin{remark}
\label{grouplattice}
Let us remark that $M_b$ and~$M_r$ are both global quotients and
that this def\/ines their quasifold structures.
The quasilattice~$F$ can be viewed as the group quasilattice of the quasitorus~$D^3$ acting on both.
\end{remark}

Remark now that, by Proposition~\ref{rotazioni}, {\em each} of the
oblate and~prolate rhombohedra in the tiling can be obtained from~$\db$ and~$\dr$ respectively by a~transformation of the icosahedral
group composed with a~translation.
We can then prove the following

\begin{theorem}
\label{uguali}
Consider an Ammann tiling having edge length $\sigma$.
Then
the compact connected symplectic quasifold corresponding to each
oblate rhombohedron in the tiling is given by~$M_b$, while the
compact connected symplectic quasifold corresponding to each prolate
rhombohedron is given by~$M_r$.
\end{theorem}

\begin{proof}
Observe that, for
each oblate rhombohedron, there exists a~transformation $T$ in the
icosahedral group that leaves the quasilattice $F$ invariant, that
sends the orthogonal vectors relative to the chosen oblate
rhombohedron to the orthogonal vectors relative to $\db$, and~such
that the dual transformation $T^*$ sends $\db$ to a~translate of the
chosen oblate rhombohedron.
The same reasoning applies to the
prolate rhombohedra of the tiling.
This implies that the reduced
space corresponding to each oblate rhombohedron of the tiling, with
the choice of orthogonal vectors and~quasilattice specif\/ied above,
is exactly $M_b$.
This yields a~unique symplectic quasifold,~$M_b$,
for all the oblate rhombohedra in the tiling.
In the same way we
prove that we obtain a~unique symplectic quasifold, $M_r$, for all
the prolate rhombohedra in the tiling.
\end{proof}
The quasifolds $M_b$ and~$M_r$ can also be constructed as complex quotients and~are
K\"ahler~\cite{cx}.

\section[Local geometry of the quasifolds $M_b$ and~$M_r$]{Local geometry of the quasifolds $\boldsymbol{M_b}$ and~$\boldsymbol{M_r}$}
\label{geoquasi}
In this section we study the equivariant geometry of the quasifolds $M_b$ and~$M_r$ in
a~neighborhood of the $D^3$-f\/ixed points.

Let us begin by describing an atlas for the quasifold $M_b$.
The
charts of this atlas are indexed by the vertices of the polytope: in
our case we f\/ind an atlas given by eight charts, each of which
corresponds to a~vertex of the oblate rhombohedron.
Consider for
example the origin: it is given by the intersection of the facets
whose orthogonal vectors are $X_1$, $X_2$ and~$X_3$.
Let $B_b$ be the ball in~$\C$
of radius $b$, namely
\begin{gather*}
B_b=\{z\in\C\,|\, |z|<b\}.
\end{gather*}
Consider the following mapping, which gives a~slice of
$\Psi_b^{-1}(0)$ transversal to the $N$-orbits
\begin{alignat*}{4}
& B_b\times B_b\times B_b && \stackrel{t_{1}}{\longrightarrow} \ &&
\big\{\vz\in\Psi^{-1}_b(0)\, |\,z_4\neq0,\, z_5\neq0,\, z_6\neq0\big\} &\\
& (z_1,z_2,z_3) && \longmapsto \ && \Big(z_1,z_2,z_3,
\sqrt{b^2-|z_1|^2},\sqrt{b^2-|z_2|^2},\sqrt{b^2-|z_3|^2}\Big). &
\end{alignat*}
This induces the homeomorphism
\begin{alignat*}{4}
& (B_b\times B_b\times B_b)/\Gamma_{1} && \stackrel{\tau_{1}}{\longrightarrow} \ && U_{1} & \\
& [\vz] && \longmapsto \ && [t_1(\vz)], &
\end{alignat*}
where the open subset $U_{1}$ of $M_b$ is the quotient
\begin{gather*}
\big\{\vz\in\Psi_b^{-1}(0)\,|\, z_4\neq0,\, z_5\neq0,\, z_6\neq0\big\}/N
\end{gather*}
and the discrete group $\Gamma_{1}$ is given by $\Gamma_{1}\simeq
N\cap(S^1\times S^1\times S^1\times\{1\}\times\{1\}\times\{1\}),$
hence
\begin{gather}
\label{gamma1}
\Gamma_{1}=\exp\left\{\left(\phi h,\phi k,\phi l\right)
\,|\, h,k,l\in\Z\right\}.
\end{gather}
The triple
$(U_{1},\tau_{1},(B_b\times B_b\times B_b)/\Gamma_{1})$ is a~chart
of $M_b$.
Analogously, we can construct seven other charts,
corresponding to the remaining vertices of the oblate rhombohedron,
each of which is characterized by a~dif\/ferent combination of the
variables.
One can show that these eight charts are compatible and
give an atlas of~$M_b$.

One can check that the moment map, locally, on the f\/irst chart is given by
\begin{gather*}
\Phi([z_1:z_2:z_3])
=\frac{\phi\alpha_5}{\langle\phi\alpha_5,U_1\rangle}|z_1|^2+\frac{\phi\alpha_6}{\langle\phi\alpha_6,U_2\rangle}|z_2|^2
+\frac{\phi\alpha_4}{\langle\phi\alpha_4,U_3\rangle}|z_3|^2
\\
\phantom{\Phi([z_1:z_2:z_3])}
=\phi\alpha_5|z_1|^2+\phi\alpha_6|z_2|^2+\phi\alpha_4|z_3|^2,
\end{gather*}
while the isotropy action of $D^3$ on $\C^3/\G_1$ is given by
\begin{alignat}{4}
& \big(D^3,\C^3/\G_1\big) && \longrightarrow \ && \C^3/\G_1 &
\nonumber\\
& ([X],(z_1,z_2,z_3)) && \longmapsto \ && \big(e^{2\pi i\phi\alpha_5(X)}z_1,
e^{2\pi i\phi\alpha_6(X)}z_2,e^{2\pi i\phi\alpha_4(X)}z_3\big). & \label{pesib}
\end{alignat}
To obtain the local expression of the moment mapping on the other seven charts it suf\/f\/ices to
replace $\phi\alpha_5$, $\phi\alpha_6$, $\phi\alpha_4$ in~\eqref{pesib} with all the possible
combinations of $\pm\phi\alpha_5$, $\pm\phi\alpha_6$, $\pm\phi\alpha_4$ respectively.
Notice that the vectors $U_1$, $U_2$, $U_3$ are three of the six generators of $F$, while
$\phi\alpha_4$, $\phi\alpha_5$, $\phi\alpha_6$ are three of the six generators of~$\frac{1}{2}I$.

An atlas for the the quasifold~$M_r$ can be constructed in the same way.
It can be shown that the moment mapping for the prolate rhomobohedron, is given, locally on the
chart corresponding to the origin, by
\begin{gather*}
\Phi([z_1 : z_2 : z_3])=
\frac{\alpha_2}{\langle\alpha_2,U_4\rangle}|z_1|^2\!+\frac{\alpha_3}{\langle\alpha_3,U_5\rangle}|z_2|^2\!
+\frac{\alpha_1}{\langle\alpha_1,U_6\rangle}|z_3|^2
=\alpha_2|z_1|^2\!+\alpha_3|z_2|^2\!+\alpha_1|z_3|^2,
\end{gather*}
while the isotropy action of $D^3$ on $\C^3/\G_1$ is given by
\begin{alignat*}{4}
& \big(D^3,\C^3/\G_1\big) && \longrightarrow \ && \C^3/\G_1 & \\
& ([X],(z_1,z_2,z_3)) && \longmapsto \ && \big(e^{2\pi i\alpha_2(X)}z_1,e^{2\pi i\alpha_3(X)}z_2,e^{2\pi i\alpha_1(X)}z_3\big). &
\end{alignat*}

Again, notice that the vectors $U_4$, $U_5$, $U_6$ are the three remaining generators of~$F$, while
$\alpha_1$, $\alpha_2$, $\alpha_3$ are the three remaining generators of $\frac{1}{2}I$.
In conclusion, the weights of the isotropy action of the quasitorus $D^3$ on a~neighboorhood of the
$D^3$-f\/ixed points for both~$M_b$ and~$M_r$ generate the quasilattice~$\frac{1}{2}I$.
Therefore~$\frac{1}{2}I$ can be thought of, in this setting, as the weight quasilattice of~$D^3$.
This is consistent with the fact that~$\frac{1}{2}I$ is dual to the group quasilattice $F$ (cf.\
Remark~\ref{duali}).
\begin{remark}
Remark that, since $\alpha_i(X)$, $\phi\alpha_i(X)$ lie in $\Z+\phi\Z$ whenever $X\in F$, and~since
the local group in each chart of~$M_b$ and~$M_r$ is equal to~$\G_1$, the above actions are well
def\/ined.
\end{remark}

\begin{remark}
If we choose as group lattice $tF$ instead, $t\in\R$, then the corresponding weight
lattice would have to be $\frac{1}{2t}I$.
But this would not be consistent with the inclusion and~projection schemes in Remark~\ref{duali}.
This is the main reason underlying our choice of the norm of the vectors~$\alpha_j$, $j=1,\ldots,6$.
\end{remark}

\section{Dif\/feotype and~symplectotype of the tiles}
\label{tipi}
The purpose of this section is to prove the following
\begin{theorem}
\label{tipo}
The quasifolds $M_b$ and~$M_r$ are diffeomorphic
but not symplectomorphic.
\end{theorem}

Before proceeding with the proof of this theorem we need a~few more facts on the local geometry of
the quasifold $M_b$.
Let us denote by $p_b$ the projection
\begin{gather*}
S^2_b\times S^2_b\times S^2_b\rightarrow M_b.
\end{gather*}
Denote by $V_{n}$ the open subset of
$S^2_{b}$ given by $S^2_b$ minus the south pole and~by $V_{s}$ the
open subset of $S^2_b$ given by $S^2_b$ minus the north pole.
Then,
on $\Psi_b^{-1}(0)$, consider the action of $S^1\times S^1\times
S^1$ given by~\eqref{esse1peresse1}.
We obtain
\begin{gather*}
V_n\times V_n\times
V_n=\big\{\vz\in\Psi_b^{-1}(0)\,|\,z_4\neq0,\,z_5\neq0,\, z_6\neq0\big\}/
\big(S^1\times S^1\times S^1\big)
\end{gather*}
and
\begin{gather*}
U_1=(V_n\times V_n\times V_n)/\Gamma.
\end{gather*}
We have the following commutative diagram:
\begin{gather}\label{diagramma}
\begin{gathered}
\xymatrix{
B_b\times B_b\times B_b\ar[r]^>>>>>>>>{t_1}\ar[d]&
\big\{\vz\in\Psi_b^{-1}(0)\,|\,z_4\neq0,\, z_5\neq0,\, z_6\neq0\big\}\ar[d]
\\
B_b\times B_b\times B_b\ar[r]^{\tilde{\tau}_1}\ar[d]_{p_{1}}&
V_n\times V_n\times V_n\ar[d]^{p_b}
\\
(B_b\times B_b\times B_b)/\Gamma_1\ar[r]^{\tau_1}&U_1.}
\end{gathered}
\end{gather}
The mapping $\tilde{\tau}_1$ is induced by the diagram and~can be
written as $\tau_n\times\tau_n\times\tau_n$, with $\tau_n \colon
B_b\rightarrow V_n$.
Observe that the mapping
\begin{alignat*}{4}
& \C && \longrightarrow \ && V_n &\\
& w && \longmapsto \ && \big[\tau_n\big(b w/\sqrt{1+|w|^2}\big)\big] &
\end{alignat*}

\newpage

\noindent
is just  the stereographic projection from the north pole.
We denote
by $\tau_s$ the analogous mapping
$\tau_s \colon B_b\longrightarrow V_s$.
The two charts
$(B_b,\tau_n,V_n)$ and~$(B_b,\tau_s,V_s)$ give a~symplectic atlas of
$S^2_b$, whose standard symplectic structure is induced by the
standard symplectic structure on~$B_b$.
Analogously, at a~local
level, the symplectic structure of the quotient~$M_b$ is induced by
the standard symplectic structure on $B_b\times B_b\times B_b$.

We have already seen that the quasifold~$M_b$ is a~global quotient
of a~product of three $2$-spheres by the discrete group~$\Gamma$.
We remark that the atlas above is the quotient by $\Gamma$ of the
atlas of the product of three spheres, given by the eight triples
$V_n\times V_n\times V_n$, $V_n\times V_n\times V_s$,
$V_n\times V_s\times V_n$, $V_s\times V_n\times V_n$,
$V_n\times V_s\times V_s$, $V_s\times V_n\times V_s$,
$V_s\times V_s\times V_n$, $V_s\times V_s\times V_s$.

We are now ready to prove Theorem~\ref{tipo}:
\begin{proof} Let us begin by showing
that $M_b$ and~$M_r$ are dif\/feomorphic.
Let us denote by $p_r$ the projection
\begin{gather*}
S^2_r\times S^2_r\times S^2_r\rightarrow M_r.
\end{gather*}
The natural
$\G$-equivariant dif\/feo\-mor\-phism $f^{\dagger} \colon S^2_b\times
S^2_b\times S^2_b\rightarrow S^2_r\times S^2_r\times S^2_r$ induces
a homeomorphism $f \colon M_b\rightarrow M_r$; in general, a
homeomorphism between two global quotients that is induced by an
equivariant dif\/feo\-mor\-phism of the manifolds turns out to be a
quasifold dif\/feo\-mor\-phism~\cite[Def\/inition A.2]{kite}.

\looseness=1
Let us now show that $M_b$ and~$M_r$ are not symplectomorphic.
Denote by $\omega_b$ and~$\omega_r$ the symplec\-tic forms of $M_b$
and $M_r$ respectively.
Suppose that there is a~symplectomorphism
$h \colon M_b\longrightarrow M_r$, namely a~dif\/feo\-mor\-phism $h$
such that $h^*(\omega_r)=\omega_b$.
We prove that this implies that
the homeomorphism $h \colon M_b\rightarrow M_r$ lifts to a
symplectomorphism $\tilde{h} \colon S^2_b\times S^2_b\times
S^2_b\rightarrow S^2_r\times S^2_r\times S^2_r$, leading thus to a
contradiction: such symplectomorphism cannot exists, since the two
manifolds have dif\/ferent symplectic volumes.
To start with recall
from~\cite[Remark~2.9]{kite} that, to each point $m\in M_b$, one can
associate the groups $\Gamma_m$ and~$\Gamma_{h(m)}$.
The def\/inition
of dif\/feo\-mor\-phism implies that these two groups are isomorphic.
Let
$n_b\in S^2_b$ be the north pole and~take $m_0=p_b(n_b\times
n_b\times n_b)$.
Then, since $\G_m\simeq\G_{h(m)}$, without loss of
genera\-li\-ty the point $h(m_0)$ can be taken to be $p_r(n_r\times
n_r\times n_r)$, where $n_r\in S^2_r$ is the north pole.
Consider
the chart $U_{1}$ that we constructed above.
Then, by def\/inition of
quasifold dif\/feo\-mor\-phism~\cite[Def\/inition~A.23]{kite} and
\cite[Remark~A.24]{kite}, there exists an open subset $U\subset U_1$
such that $m_0\in U$ and~$h\circ\tau_1^{-1} \colon
\tau_1^{-1}(U)\rightarrow h(U)$ is a~dif\/feo\-mor\-phism of the universal
cove\-ring models induced by $\tau_1^{-1}(U)\subset B_b\times
B_b\times B_b/\G_1$ and~$h(U)\subset M_r$ respectively.
Moreover, by
\cite[Proposition~A.9]{kite}, any open subset $W\subset U$ enjoys
the same property.
We can choose $W_0\subset U_1$ such that
$\tilde{W}_{0}=(\tau_{1}\circ p_{1})^{-1}(W_0))$ is a~product of
three balls.
In particular, $\tilde{W}_{0}$ is simply connected.
Denote now by $\wt_{r,0}=(p_r)^{-1}(h(W)))$; this is an open subset
of $S^2_r\times S^2_r\times S^2_r$, which is also connected, due to
the action of $\G$ on $S^2_r\times S^2_r\times S^2_r$.
Denote by
$\wsh_{r,0}$ its universal covering.
Now consider a~point $\vz^1\in
B_b\times B_b\times B_b$ such that $z_1^1\neq0$, $z_2^1\neq0$ and
$z_3^1\neq0$ and~let $m=(\tau_1\circ p_1)(\vz^1)$.
For the sequel it
is crucial to remark that, because of the action of $\G_{1}$ given
in~\eqref{gamma1}, any $\G_{1}$-invariant open subset of $B_b\times
B_b\times B_b$ that contains the point $\vz^1$, contains also the
product of circles $\{(z_1,z_2,z_3)\in B_b\times B_b\times
B_b\,|\,|z_1|=|z^1_1|,\,|z_2|=|z^2_1|,\,|z_3|=|z^3_1|\}$.
Hence, for
each point $(\tau_{1}\circ p_1)(t\vz^1)$ with $t\in[0,1]$, we can
f\/ind an open subset $W_t\subset U_1$, containing that point, such
that the homeomorphism $\tau_1^{-1}\circ h$, restricted to
$\tau_1^{-1}(W_t)$, is a~dif\/feo\-mor\-phism, and~$(\tau_1\circ
p_1)^{-1}(W_t)$ is the product of three open annuli.
We can cover
the curve by a~f\/inite number of these $W_t$'s: $W_0,W_1,\ldots,W_s$,
with $W_j\cap W_{j+1}\neq\varnothing$.
Notice that $(\tau_1\circ
p_1)^{-1}(W_j\cap W_{j+1})$, $j=0,\ldots,s-1$, is itself a~product
of three open annuli.
The subsets $\wt_{j}=(\tau_1\circ
p_1)^{-1}(W_j)$ and~$\wt_{r,j}=(p_r)^{-1}(h(W_j))$ are open and
connected.

We divide the remaining part of the proof in subsequent steps:

{\it Step 1}: consider f\/irst $W_0$.
Since the isotropy of
$\G_1$ at $0$ is the whole $\G_1$, we can apply
\cite[Lemma~6.2]{kite}.
We f\/ind that $\wt_{r,0}$ is itself simply
connected and~that the homeomorphism $h\circ\tau_1$ lifts to a
dif\/feo\-mor\-phism $\tilde{h_0}\,\colon\,\wt_{0}\rightarrow\wt_{r,0}$.

{\it Step 2}: consider the homeomorphism $h_1=h\circ\tau_{1}$
def\/ined on $\tau_{1}^{-1}(W_1)$.
By construction $h_1$ is a
dif\/feo\-mor\-phism of the universal covering models of the induced
models.
We f\/ind the following diagram:
\begin{gather*}
\xymatrix{
\wsh_{1}\ar[r]^{h^{\sharp}_1}\ar[d]_{\pi_1}&\wsh_{r,1}\ar[d]^{\rho_1}
\\
\wt_{1}\ar[d]_{p_1}&\wt_{r,1}\ar[d]^{q_1}
\\
\tau_{1}^{-1}(W_1)\ar[r]^{h_1}&h(W_1).
}
\end{gather*}
Consider the restriction of $h_1$ to $\tau_{1}^{-1}(W_0\cap W_1)$.
This restriction admits a~lift, given by the restriction of
$h^{\sharp}_1$ to $(\pi_1\circ p_{1})^{-1}(\tau_{1}^{-1}(W_0\cap
W_1))$.
Furthermore, by Step~1, the restriction of $h_1$ admits
another lift, def\/ined on $p_{1}^{-1}(\tau_{1}^{-1}(W_0\cap W_1))$,
which is the restriction of $\tilde{h}_0$.
Therefore, by
\cite[Lemma~6.3]{kite}, the restriction of $\rho_1\circ
h_1^{\sharp}$ to $(\pi_1\circ p_{1})^{-1}(\tau_{1}^{-1}(W_0\cap
W_1))$ descends to a~dif\/feo\-mor\-phism def\/ined on
$p_{1}^{-1}(\tau_{1}^{-1}(W_0\cap W_1))$.

{\it Step 3}: we consider $W_0\cap W_1\subset W_1$ and~we
apply~\cite[Lemma~6.5]{kite} to the homeomorphism $h\circ\tau_{1}$
def\/ined on $\tau_{1}^{-1}(W_1)$.
We deduce that $h\circ\tau_{1}$ is
a dif\/feo\-mor\-phism of the model $(\tau_{1}\circ
p_{1})^{-1}(W_1)/\G_{1}$
with the model induced by $h(W_1)\subset M_r$.

{\it Step 4}: we apply Step~3 to the other successive
intersections.
We f\/ind that $h\circ\tau_{1}$ is a~dif\/feo\-mor\-phism of
the model $(\tau_{1}\circ p_{1})^{-1}(\mathop{\cup}\limits_{i=1}^k W_i)/\G_{1}$ with
the model induced by $h(\mathop{\cup}\limits_{i=1}^k W_i)\subset M_r$.
Remark now that a~slight modif\/ication of the above argument applies to any
choice of point
$\vz^1\in B_b\times B_b\times B_b$, $\vz_1\neq0$.

Let $\epsilon>0$ be arbitrarily small.
Consider the product of
closed balls $\overline{B}_{b-\epsilon}\times
\overline{B}_{b-\epsilon}\times\overline{B}_{b-\epsilon}$.
This, by
Step~4, can be covered by a~f\/inite number of connected open subsets
of the kind $(\tau_{1}\circ p_{1})^{-1}(\mathop{\cup}\limits_{i=1}^k W_i)/\G_{1}$,
whose intersection is a~product of three balls centered at the
origin.
Now~\cite[Lemma~A.3]{kite}, which guarantees the uniqueness
of the lift up to the action of $\G$, implies that the homeomorphism
$h$ admits a~lift to $\tilde{\tau}_1(\overline{B}_{b-\epsilon}
\times\overline{B}_{b-\epsilon}\times\overline{B}_{b-\epsilon})$.
This in turn implies that $h \colon U_1\rightarrow h(U_1)$ admits
a lift
\begin{gather*}
\tilde{h}_1 \colon \  V_n\times V_n\times V_n\rightarrow p_r^{-1}(h(U_1)).
\end{gather*}
We apply the same argument to the other eight charts.
These charts
intersect on the dense connected open subset where the action of the
quasitorus~$D^3$ is free.
By the uniqueness of the lift \cite[Lemma~A.3]{kite}, we obtain a~global lift
$\tilde{h}\,\colon\,S^2_b\times S^2_b\times S^2_b\rightarrow
S^2_r\times S^2_r\times S^2_r$.
Moreover, since diagram
\eqref{diagramma} preserves the symplectic structures, we have that
$\tilde{h}$ is a~symplectomorphism between $S^2_b\times S^2_b\times
S^2_b$ to $S^2_r\times S^2_r\times S^2_r$, which is impossible.
\end{proof}

In conclusion, similarly to what happens in dimension two for Penrose rhombus
tilings~\cite{rhombus}, there is a~{\em unique} quasifold structure that is
naturally associated to any Ammann tiling with f\/ixed edge length,
and two distinct symplectic structures that distinguish the oblate
and the prolate rhombohedra.

\subsection*{Acknowledgements}
We would like to thank Ron Lifshitz for his help on the theory of quasicrystals.

\pdfbookmark[1]{References}{ref}
\LastPageEnding

\end{document}